\documentclass{amsart}
\usepackage[all]{xy}
\usepackage{graphicx}
\usepackage{psfrag}
\usepackage{amsthm}
\usepackage{amscd}
\usepackage{amsfonts}
\usepackage{amstext}
\usepackage{amssymb}
\usepackage{amsmath}
\usepackage{amsxtra}
\usepackage{plain}
\usepackage[all]{xy}

\usepackage[active]{srcltx}
\usepackage{colordvi}

\newtheorem{prop}{\bf Proposition}[section]

\newtheorem{lem}[prop]{{\bf Lemma}}
\newtheorem{thm}[prop]{{\bf Theorem}}
\newtheorem{rem}[prop]{{\bf Remark}}

\numberwithin{equation}{section}

\newenvironment{pf}{{\bf proof }}{\qed\endtrivlist}

\newcommand{\Z}{\mathbb{Z} }

\newcommand{\R}{\mathbb{R} }

\newcommand{\tr}{\operatorname{tr}}

\newcommand{\ind}{\operatorname{ind}}

\newcommand{\im}{\operatorname{Im}}

\newcommand{\sig}{\operatorname{sig}}

\long\def\symbolfootnote[#1]#2{\begingroup%
\def\thefootnote{\fnsymbol{footnote}}\footnote[#1]{#2}\endgroup}

\begin{document}

\title{On cut-and-Past invariance of Kervaire semi-characteristic}
\author{Mostafa ESFAHANI ZADEH}

\address{Mostafa Esfahani Zadeh.\newline
Mathematical Sciences Department, 
Sharif University of Technology,\newline
Tehran-Iran\and\newline
School of Mathematics, Institute for research in Fundamental Sciences (IPM)\newline
P. O. Box:19395-5746, Tehran, Iran}

 \email{esfahani@sharif.ir}

\symbolfootnote[0]{\emph{2000 Mathematics Subject Classification}. 58E05, 35J25.} 
\symbolfootnote[0]{\textit{Key words and phrases}. Witten deformation, Kervaire invariant, Index Theory, 
Cut-and-Past operation.}
\symbolfootnote[0]{This research was in part supported by a grant from IPM (NO. 89510130)}

\begin{abstract}
In this note we study the relative Kervaire semi-characteristic and prove its 
invariance under cut-and-past operation. Our approach is analytic and follow very closely 
the method introduced by W. Zhang
\end{abstract}
\maketitle

\section{Introduction}
The purpose of this paper is to prove the Cut-and past invariance of Kervaire semi-characteristic for 
orientable manifolds of dimension $n=4q+1$. The method we use is analytic and is introduced by G. Zhang 
in \cite{Zhang2}. In fact the main step is generalizing the theorem 1.3 of the previous 
reference to manifolds with boundary. We 
describe here briefly the content of the paper \cite{Zhang2}. 
Let $M$ be a $n$-dimensional \emph{closed}  manifold with
vanishing Euler character, $\chi(M)=0$. This last condition implies the existence of a nowhere-vanishing 
Vector field $V$ on $M$. 
Let $E$ be a vector bundle over $M$ whose fiber at $x$ is a complementary subspace for 
$V(x)\in T_xM$. If $X$ is a transversal section of $E$ then the set $Z$ of its singularities is 
a one-dimensional closed sub-manifold of $M$. Therefore $Z$ is a disjoint union of some 
closed curves $\gamma$'s. 
The vector field $X$ is used to construct a real line bundle $\mathcal L$ over each closed curve $\gamma$ 
and the index $\ind_2(\gamma)\in Z_2$ of the curve $\gamma$ is defined as follows  
\[\ind_2(\gamma)=\left\{
\begin{array}{rl}
1&\text{ if }\mathcal L\text{ is trivial}\\
0&\text{ if }\mathcal L\text{ otherwise}
\end{array}\right.
\] 

The main results of \cite{Zhang2} assert that the sum $\sum_{\gamma}\ind_2(\gamma)$ 
is a topological invariant and the following equalities hold in $\Z_2$
\begin{equation}\label{zhangasli}
\sum_{\gamma}\ind_2(\gamma)=\left\{
\begin{array}{ll}
0&\text{ if }\dim M=4q+2\text{ or }4q+3\\
\sig (M)/2&\text{ if }\dim M=4q\\
\kappa(M)&\text{ if }\dim M=4q+1
\end{array}\right.
\end{equation}
Here for an oriented $4q$-dimensional closed manifold $\sig(M)$ denotes the signature of $M$ while   
$\kappa(M)\in\Z_2$ stands for the Kervaire semi-characteristic of $M$ which is defined by the 
following relation 
\[\kappa(M)=\sum_j\dim H^{2j}(M)\in \Z_2.\] 
In above as in the body of the paper we deal with the de Rham cohomology spaces. The reader 
should notice the analogy of the above equality with the Poincaré-Hopf formula for the Euler 
characteristic. Since this formula holds only where the Euler Characteristic vanishes, it might be 
considered as a secondary Poincaré-Hopf formula. The signature of the $4q$-dimensional manifold is 
by definition the signature of the following symmetric non degenerate bilinear form
\begin{gather*}
\beta:H^{2q}(M)\times H^{2q}(M)\to\R\\
\beta(\omega,\eta)=\int_M\omega\wedge\eta~.
\end{gather*}
A remarkable property of the signature $\sig (M)$ is expressed by the Hirzebruch 
signature theorem. This theorem asserts that $\sig (M)$ is a characteristic number and, therefore, 
a cobordism invariant. In particular if $M$ is the boundary of another manifold then 
$\sig (M)=0$. 

In this paper we prove that the Kervaire semi-characteristic of a closed orientable 
$(4q+1)$-dimensional manifold is a cut-and past invariant. More precisely let $M$ be a closed manifold without 
boundary and let $N\subset M$ 
be a sub-manifold of co-dimension one that partitions $M$ into two compact manifolds $M_1$ and $M_2$ 
with common boundary $\partial M_1=\partial M_2=N$. Given a diffeomorphism  $\phi: N \to N$ 
the closed manifold 
$M_1\sqcup_\phi M_2$ is defined by identifying points $x\in \partial M_1$ and $\phi(x)\in\partial M_2$. 
We prove in theorem \ref{theorem2} that in dimension $n=4q+1$ the equality 
$\kappa(M)=\kappa(M_1\sqcup_\phi M_2)$ holds provided that $\chi(N)=0$. To prove this 
theorem we state and prove a generalization of the counting Zhang formula for manifolds with boundary. 

Let $M$ be a compact manifold of dimension $n=4q+1$ with non-empty boundary $\partial M$. We assume that 
$\chi(\partial M)=0$ which implies the existence of a nowhere vanishing vector field $V$ on $\partial M$. 
Since the dimension of $M$ is odd, this vector field has an extension to a nowhere vanishing vector 
filed on $M$ that we denote by the same symbol $V$. We define the bundle $E$ as in the boundary-less case.  
Let $X$ be a transverse section of $E$ which is tangent to the boundary. 
The set of all singularities of $X$ is the union of a finite number of 
closed curves $\gamma$'s. This closed curves are either subsets of $M^\circ$; the interior of $M$; 
or subsets of the boundary.  
We will construct the real line bundles $\mathcal L$ over $\gamma$'s. 
The index $\ind_2(\gamma)\in\Z_2$ is defined as before. 
Define the relative Kervaire semi-characteristic $\kappa(M,\partial M)$ by the following relation 
\begin{equation}\label{kerver}
 \kappa(M,\partial M)=\sum_{j}\dim H^{2j}(M,\partial M)\in\Z_2.
\end{equation}
We prove the following generalization of the Zhang formula in this context 
\[\kappa(M,\partial M)=\sum_{\gamma\subset M^\circ}\ind(\gamma).\] 
As it is clear from this formula the singular curves in boundary have no contribution in calculating the 
relative Kervaire invariant.
   
{\bf Acknowledgement: } The author would like to thank Mathematics Institute of the 
Georg-August university in Goettingen. A part of this paper is written when the author 
was a visitor there.  

\section{The geometric and analytic setting }

Let $M$ be a compact manifold with boundary $\partial M\neq\emptyset$. Let $\Omega^p(M,\partial M)$ 
denote the set of differential $p$-forms $\omega$ such that $j^*\omega\in \Omega^p(\partial M)$ 
is exact. It is clear by this definition that the exterior derivative $d$ maps 
$\Omega^p(M,\partial M)$ into $\Omega^{p+1}(M,\partial M)$. Therefore one has the 
following differential complex 
\[\dots\to\Omega^{p-1}(M,\partial M)\stackrel{d}{\rightarrow}\Omega^p(M,\partial M)
\stackrel{d}{\rightarrow}
\Omega^{p+1}(M,\partial M)\to\dots\]
which define the relative de Rham cohomology groups 
\[H^p(M,\partial M):=\dfrac{\ker d:\Omega^p(M,\partial M)\to\Omega^{p+1}(M,\partial M)}
{\im d:\Omega^{p-1}(M,\partial M)\to\Omega^{p}(M,\partial M)}.\]
For giving a Hodge theoretic interpretation of these groups it 
is necessary to have a Riemannian metric $g$ on $M$. We assume that this metric takes 
the following product form in a collar neighborhood $U=(-1,0]\times\partial M$
\begin{equation}\label{prod}
g(u,y)=d^2u+g_0(y),
\end{equation}
where $g_0$ is a Riemannian metric on $\partial M$. This Riemannian structure induces an 
inner product on $\Omega^*(M)$. The symbol $\delta:\Omega^*(M)\to\Omega^{p-1}(M)$ 
denotes the formal adjoint of the operator $d$ with respect to this inner product. 
Each differential form $\omega$ takes the following form in the collar neighborhood $U$
\[\omega_{|U}=\omega_0(u,y)+du\,\wedge\omega_1(u,y).\]
The differential form $\omega$ satisfies the relative boundary condition if the following 
relations hold 
\begin{equation}\label{diri}
\omega_0(0,y)=0~~\text{ and }~~\frac{\partial\omega_1}{\partial t}(0,y)=0
\end{equation}
Let $\Omega^*(M,B)\subset\Omega^*(M)$ consist of all differential forms satisfying 
the relative boundary condition.  
For $\eta=\eta_0+du\,\wedge\eta_1$ another differential form  and $D=d+\delta$ one has the 
Green formula, c.f. \cite[page 24]{BoWo1}
\begin{equation}\label{green}
\langle D\omega,\eta\rangle-\langle\omega,D\eta\rangle
=\int_{\partial M}\langle\omega_0,\eta_1\rangle-\int_{\partial M}\langle\omega_1, \eta_0\rangle,
\end{equation}
This implies that the operator $D$ is formally self-adjoint on $\Omega^*(M,B)$. 
Therefore the Laplacian operator $\triangle =D^2$ is a formally positive second order 
elliptic operator on $\Omega^*(M,B)$ and its kernel consists of smooth closed differential forms satisfying 
the conditions \eqref{diri} and  
\[\ker\triangle=\ker D.\]
The following relation is the consequence of the Hodge theory in this relative boundary 
condition context 
\begin{equation}\label{hodgrel}
H^*(M,\partial M)\sim\ker \triangle 
\end{equation}
The spaces $\Omega^p(M,B)$; for $0\leq p\leq n$; can be completed with respect to appropriate Sobolev 
norm to obtain the Sobolev spaces $W^p(M,B)$. The inclusion $\Omega^*(M,B)\hookrightarrow\Omega^*(M)$
induces an isometric embedding $W^p(M,B)\hookrightarrow W^p(M)$. Since the image of this embedding is closed,   
all classical theorems in the theory of the Sobolev 
spaces, e.g. Sobolev embedding theorem, Rellich's theorem and the elliptic estimate 
hold also in this context.

For next uses we need to introduce some tools from linear and Clifford algebra. For this purpose 
let $E$ be a real vector space with an inner product $\langle,\rangle$. For $v\in E$ 
consider the following linear maps on $\Lambda^*E^*$  
\[c(v):=v^*\wedge-i_v~~;~~\hat c(v)=v^*\wedge+i_v,\]
where $v^*\in E^*$ denotes the dual of the vector $v$ with respect to the inner product. 
The inner product of $E$ induces a natural inner product on $\Lambda^*E^*$ and these 
operators satisfy the following properties with respect to this inner product 
\begin{equation}\label{triup}
c(v)^*=-c(v)~~;~~\hat c(v)^*=\hat c(v)^*.
\end{equation} 
Moreover it is easy to verify the following relations 
\begin{align}\label{clirels}
c(v)c(w)+c(w)c(v)&=-2\langle v,w\rangle,\notag\\
\hat c(v)\hat c(w)+\hat c(w)\hat c(v)&=2\langle v,w\rangle,\\
c(v)\hat c(w)+\hat c(w)c(v)&=0\notag
\end{align}
Let $A:E\to E$ be an invertible self adjoint linear map whose matrix with respect to 
an orthonormal basis $\lbrace e_i\rbrace$ is $a_{ij}$. 
The following endomorphism on $\Lambda^*E^*$ is independent of the basis  
\[\hat A=\sum_{ij}a_{ij}c(e_i)\hat c(e_j).\]
The following lemma gives a property of this operator that will  be important in the next section. 
\begin{lem}\cite[lemma 4.8]{Zhang1}\label{lem1}
The following operator 
\begin{gather*}
K:\Lambda^*E^*\to \Lambda^*E^*\\
K=\tr|A|+\hat A
\end{gather*}
is non-negative, its kernel is one dimensional and is generated by an element in $\Lambda^{even}E^*$ 
if $det A>0$ and in $\Lambda^{odd}E^*$ if $\det A<0$. 
\end{lem}
In fact because $A$ is symmetric it can be diagonalized in an orthonormal basis $\{e^j\}$. If the corresponding 
eigenvalues are positive for $j=1,\dots p$ and negative for $j=p+1\dots n$ then the kernel of $K$ is 
generated by $e^{p+1}\wedge e^{p+2}\wedge\dots\wedge e^{n}$.

\section{Counting formula}

Before starting the main discussion we recall a fact from differential topology. It is true that if the 
Euler character of a manifold (with or without boundary) vanishes then this manifold carries a nowhere 
vanishing vector field which is outward on boundary points, 
c.f. \cite[chapter 6]{Mi1}. Let
$(M,\partial M)$ be an odd dimensional manifold with boundary. The Poincaré duality implies 
the relation $\chi(M,\partial M)=-\chi(M)$. On the other hand the long exact sequence of the pair 
$(M,\partial M)$ in de Rham context reads
\[\dots\to H^j(M,\partial M)\to H^j(M)\to H^j(\partial M)\to H^{j+1}(M,\partial M)\to\dots\]   
This implies, after summation over alternate traces, the relation 
\[\chi(M,\partial M)-\chi(M)+\chi(\partial M)=0.\] 
From now-on we assume the following condition on the boundary manifold 
\begin{equation}\label{conbon}
\chi(\partial M)=0.
\end{equation} 
This condition along the above discussion implies the vanishing of $\chi(M,\partial M)$ and $\chi(M)$. 
Therefore there is a nowhere vanishing vector field $V'$ on $M$ which is outward at boundary points. 
On the other hand, and with condition \eqref{conbon}, there is a nowhere vanishing vector field $V_0$ 
on $\partial M$. We can extend this vector field in a collar neighborhood $U=(-1,0)\times \partial M$ 
of $\partial M$ such that its restriction to $\lbrace t\rbrace\times\partial M$ is $V_0$ with respect to 
the obvious identification $(t,x)\sim x$. We also denote this extension by $V_0$. If $\eta$ is a smooth 
function on $M$ which coincides to $t$ near boundary and is supported in a small neighborhood 
of the boundary then the  vector field 
$V=\eta(t)V'+(1-\eta(t))V_0 $
is a nowhere vanishing vector field on $M$ which is tangent to the boundary.

Let $M$ be equipped with a Riemannian metric $g$ which takes the product form \eqref{prod} in the 
collar neighborhood $U$ of the boundary. From now on we assume the condition 
\eqref{conbon} is satisfied, therefore there is a nowhere vanishing vector filed $V$ on $M$ such that 
$V(x)\in T_x\partial M$ for $x\in\partial M$. We can and will assume that $\|V\|=1$. 
For each $x\in M$ let $E_x\subset T_xM$ denote the 
orthogonal compliment of $V(x)\subset T_xM$ and denote the set of all $E_x$'s 
by $E$ which is a vector bundle of rank $4q$ over $M$. 
Let $X$ be a transverse section of $E$. In the boundary points we assume $X_{|\partial M}$ be 
tangent to the boundary. 
Let $X=(X_1,X_2,\dots, X_n)$ with respect to collar coordinates \eqref{prod}. 
With above conditions in boundary points  $\partial X_1/\partial u\neq0$. 
Therefore we may assume that 
\begin{equation}\label{vecboun}
\frac{\partial X_1}{\partial u}=\pm 1 
\end{equation}
The zero set $Z$ of $X$ is a compact $1$-dimensional 
sub-manifold of $M$ whose connected components, i.e. topological circles, are either 
disjoint from or included in the boundary. In the sequel we denote these circle generically by $\gamma$.

By transversality of $X$, 
the vector field $V$ is no-where perpendicular to $\gamma$. 
In fact $\ker T_xX=T_x\gamma$ and by transversality $T_xM=\ker T_x X\oplus E_x$ 
for $x\in\gamma$, so $T_xM=T_x \gamma\oplus E_x$.  This implies, since $V$ is 
everywhere orthogonal to $E$, that $V$ is no-where orthogonal to $\gamma$. Consequently one 
can deform continuously the vector field $V$ to a unit vector field which is tangent to  
curves  $\gamma$'s. We denote this new vector field by $V$  and 
deform the Riemannian metric in neighborhoods of curves so that 
the vector field $V$ be orthogonal to $E$. 
We work under these assumptions because the quantities that we are interested 
in are homotopy invariant. \newline 
Let $x$ be a singular point of $X$. By transversality of $X$, the map $T_xX:E_x\to E_x$ 
is a linear isomorphism that we denote by $A(x)$. If $x$ is a boundary point then $T_xX:E_{0x}\to E_{0x}$ 
is also an isomorphism and we denote this isomorphism by $A_0(x)$. Let $\lbrace e_j\rbrace_{j=2}^{n}$ 
be an orthonormal basis of 
$E_x$ and $A_{ij}(x)$ be the matrix of $A(x)$ with respect to this basis. The following 
operator is indeed basis independent
\begin{gather*}
\hat A(x):\Lambda^*E_x^*\to\Lambda^*E_x^*\\
\hat A(x)=\sum_{ij=2}^{n}A_{ji}(x)c(e_i)\hat c(e_j). 
\end{gather*}
If $x$ is a boundary points then we assume that $e_2=\frac{\partial}{\partial u}$. In this case 
the following operator is basis independent 

\begin{gather*}
\hat A_0(x):\Lambda^*E_{0x}^*\to\Lambda^*E_{0x}^*\\
\hat A_0(x)=\sum_{ij=3}^{n}A_{0ji}(x)c(e_i)\hat c(e_j). 
\end{gather*}

Therefore, given a circle $\gamma$, have a field $A$ (or $A_0$ if $\gamma$ is in the boundary) of 
invertible matrices of order $(n-1)\times(n-1)$ (or of order $(n-2)\times(n-2)$) along $\gamma$  as well as 
a field of Clifford morphisms $\hat A$ (or $\hat A_0$) acting on sections of the bundle 
$\Lambda^*E^*$ (or $\Lambda^*E_0^*$) over $\gamma$. 
From Lemma \ref{lem1} for each $x\in \gamma$ the operators  
\begin{gather}\label{defk}
K(x):=\tr |A(x)|+\hat A(x): \Lambda^*E_x^*\to\Lambda^*E_x^*\\
K_0(x):=\tr |A_0(x)|+\hat A_0(x): \Lambda^*E_{0x}^*\to\Lambda^*E_{0x}^*
\end{gather}
are non-negative with one-dimensional kernels. Moreover, by the discussion after lemma \ref{lem1} and special 
form of the vector field $X$ in the collar neighborhood $U$ (which is expressed by \eqref{vecboun}), 
for $x\in \partial M$  one of the relations 
$\ker K_0(x)=j^*\ker K(x)$ or $\ker K_0(x)=i_{e_2}\ker K(x)$
hold depending on the sign in \eqref{vecboun} (here $j:\partial M\to M$ is the inclusion map). 
According to above discussion the kernel of the operators $K$ (and $K_0$) form a  
real line bundle $\mathcal L$ over $\gamma$. 
In \cite{Zhang2} the index $\ind_2(\gamma)$ is defined to be $1$ if 
this bundle is trivial and $0$ if not. Moreover the relations \eqref{zhangasli} are 
proved. 

\begin{rem}\label{asalio}
If $\gamma\subset\partial M$ then the line bundles made of $\ker K$ and $\ker K_0$ 
are naturally isomorphic. Therefore the index $\ind_2(\gamma)$ is the same if we define it 
through the boundary data. 
\end{rem}

In particular for orientable closed manifolds of dimension $n=4q+1$ the Kervaire semi-characteristic 
is proved to be given by the sum of the indices of the singular circles $\gamma$'s. 
We are going to prove 
the similar theorem for the case where the boundary is not empty.

\begin{thm}\label{aslio}
Let $M$ be a $n=4q+1$ dimensional orientable manifold with boundary $\partial M$ which satisfies 
the relation $\chi(\partial M)=0$. The following formula holds in $\Z_2$ 
\[\kappa(M,\partial M)=\sum_{\gamma\in M^{\circ}} \ind_2(\gamma)\]
Here the summation is taken over those circles which are included 
in $M^\circ$; the interior of $M$; and $\kappa(M,\partial M)$ is given by \eqref{kerver}
\end{thm}

In particular the boundary has no contribution in calculating the relative 
Kervaire semi-characteristic. The rest of this section is devoted to prove this 
theorem. We follow very closely 
the method introduced by Zhang in \cite{Zhang2}. 
We recall that $\Omega^*(M,B)$ consists of the differential forms on $M$ which satisfy the relative  
boundary condition \eqref{diri}. Since the vector field $V$ is tangent to the boundary, the Clifford actions 
$c(V)$ and $\hat c(V)$ on $\Omega(M)$ respect this boundary condition and define actions on 
$\Omega(M,B)$. Therefore the following operator acts on $\Omega^{even}(M,B)$
\[D_V=\frac{1}{2}(\hat c(V)(d+\delta)-(d+\delta)\hat c(V))\]
This is indeed a first order elliptic differential operator which is, thank to relations \eqref{green} and 
\eqref{triup}, skew-adjoint. So its Kernel is finite dimensional and its dimension modulo 
$2$ is invariant with respect 
to continuous deformation of $D_V$, c.f.  \cite{AtSiskew}. This invariance is denoted by $\ind_2(D_V)$.
The relevance of this operator is through the following formula 
\[\ind_2(D_V)=\kappa(M,\partial M).\]
The proof is completely similar to the proof of theorem 1.2 in \cite{Zhang2} once one the 
relation  \eqref{hodgrel}.

Let $s$ be a parameter varying in the set of non-negative real numbers. The vector filed $X$ can be used 
to define the following deformation of the operator $D_V$
\begin{equation}\label{defnah}
D_s=D_V+s\hat c(V)\hat{c}(X):\Omega^{even}(M,B)\to\Omega^{even}(M,B)
\end{equation}
Because $V$ and $X$ are orthogonal, It is clear from relations \eqref{clirels} that 
$\hat c(V)\hat{c}(X)$ is a skew-adjoint operator on $\Omega^*(M,B)$ of order zero. 
Therefore $\ind_2(D_s)$ is well defined and independent of $s$, i.e. 

\begin{equation}\label{indin}
\ind_2(D_s)=\ind_2(D_V)=\kappa(M,\partial M)
\end{equation}
Let $\lbrace e_j\rbrace$ for $1\leq j\leq n$ be a local oriented orthonormal basis for $TM$ and denote 
by $\nabla$ the Riemannian connection on $TM$, as well as its extension to $\Lambda^*T^*M$.
The following formula is a direct result of relations \eqref{clirels}
\begin{equation}\label{Bochner}
D_s^2=D_V^2-s\sum_{j=0}^{4q}(c(e_j)\hat c(\nabla_{e_j}X)-\langle\nabla_{e_j}X,V\rangle c(e_j)\hat c(V))-s^2|X|
\end{equation}
Let $Z$ be the set of the singular points of $X$ and $W$ an open neighborhood of $Z$. The following lemma 
can be proved using the classical theorem of Sobolev spaces. For a proof in boundary-less case we 
refer to \cite[proposition]{Zhang1}. This proof goes on literary in our boundary case.  

\begin{lem}
There exist positive constant $a$ and $b$ such that for any $s\geq 1$ and $\omega\in\Omega^*(M,B)$ 
the following 
inequality holds 
\[\|D_s\|_0^2\geq a(\|\omega\|_1^2+(s-b)\|\omega\|_0^2\]
provided that $\text{supp}(\omega)\cap W=\emptyset$. 
\end{lem}

This lemma suggests that to compute the dimension of $\ker D_s$, when $s$ gets larger and larger, it suffices 
to compute the dimension of the kernel of the restriction  of $D_s$ to arbitrary small open 
neighborhoods of $Z$. Even more, to compute the dimension of the kernel of this operator when $s$ 
goes toward infinity, 
it is enough to compute the dimension of the kernel of the corresponding linearized operator along 
closed curves $\gamma$'s. 
Fix a circle $\gamma\subset Z$. Since $M$ is oriented the normal bundle $E\to \gamma$ is oriented and trivial. 
so $E\simeq \gamma\times \R^{4q}$.  
Let $(t,y_2,y_3,\dots,y_{n})$ denote the Euclidean coordinates of $\gamma\times \R^{n}$. If $\gamma$ 
is included in the boundary we let $y_2=u$ be the collar parameter of \ref{prod}. 
Using these coordinates we have the following Bochner type formula \eqref{Bochner} where $\mathrm D_s$ 
stands for the linearization of $D_s$.  
\begin{equation*}
\mathrm D_s^2=-\frac{\partial^2}{\partial t^2}+
\lbrace -\sum_{i=2}^n\frac{\partial^2}{\partial y_i^2}+
s\hat A(t)+s^2\langle|A|y,|A|y \rangle\rbrace.
\end{equation*}  

Using \eqref{defk} the above expression can be rewritten as follows 

\begin{equation}\label{intsim}
\mathrm D_s^2=-\frac{\partial^2}{\partial t^2}+\lbrace-\sum_{i=2}^n\frac{\partial^2}{\partial y_i^2}
-s\tr |A(t)|+s^2\langle|A|y,|A|y \rangle\rbrace +sK(t). 
\end{equation}

In above expression the operator in the braces is harmonic oscillator 
operator whose kernels is one dimensional generated by the function 
$\exp(-\frac{s}{2}\langle|A|y,|A|y \rangle)$. 
The operator $K(t)$ is an algebraic operator and by lemma \ref{lem1} its 
kernel is one-dimensional real subspace of $\Lambda^{even} E^*_t$ if $\det A(t)>0$ 
and of $\Lambda^{odd}E^*_t$ if $\det A(t)<0$. 
In each case these subspaces put together form a real line bundle $\mathcal L$ over $\gamma$. 
Now we assume the sign positive in \eqref{vecboun} is occurred (the negative sign is similar).  
In this case the kernel of the 
operator \eqref{intsim} has the following form if $\det A>0$ 
\begin{equation}\label{sher1}
\exp(-\frac{s}{2}\langle|A|y,|A|y \rangle)\omega,
\end{equation}
while it takes the following form if $\det A<0$ 
\begin{equation}\label{sher2}
\exp(-\frac{s}{2}\langle|A|y,|A|y \rangle)dt\wedge\omega,
\end{equation} 

$\omega$ in \eqref{sher1} is a parallel section of $\mathcal L$ of even degree while in \eqref{sher2} 
its degree is odd (this is why we have put $dt$ in the later case). So, independent of the sign of 
$\det A$, as long as the line bundle $\mathcal{L}$ is trivial and $\gamma$ is included in the 
interior of $M$, the kernel of $\mathrm D_s$, as an operator on sections of $E\to\gamma$ is one-dimensional. 
On the other hand, in the boundary points the operator $\mathrm D_s^2$ takes the following form 
\[\mathrm D_s^2=-\frac{\partial^2}{\partial t^2}-\frac{\partial^2}{\partial u^2}-\lbrace
\sum_{i=3}^n\frac{\partial^2}{\partial y_i^2}
-s\tr |A(t)|+s^2\langle|A|y,|A|y \rangle\rbrace +sK(t)\]
If $\gamma$ is a boundary curve then by our assumption on the sign in \eqref{vecboun} the differential 
form $\omega$ in \eqref{sher1} does not contain the term $du$. 
Therefore $dt\wedge du\wedge \omega$ satisfies the boundary condition \eqref{diri} and generates 
the kernel of $\mathrm D_s$. 
If $det A<0$ then $\omega$ in relation \eqref{sher2} is a parallel section of $E_\to\gamma$ with odd degree. 
In this case  $du\wedge \omega$ satisfies the boundary condition \eqref{diri}  
and generates the kernel of $\mathrm D_s$ provided that $\mathcal L$ is trivial. Summarizing the 
dimension of $\ker \mathrm D_s$ equals the number of closed curves $\gamma$'s whose associated line bundles 
$\mathcal L$ are trivial, regardless this curve is in the interior or on the boundary of $M$, i.e. 
\[\dim \ker D_s=\sum_{\gamma\subset M^\circ}\ind_2(\gamma)+\sum_{\gamma\subset \partial M}\ind_2(\gamma)\]
By remark \ref{asalio} for evaluating the boundary contribution in this formula we can just use 
the boundary date. By the formula \eqref{zhangasli} the boundary contribution in above formula is equal to 
$\sig(\partial M)/2$ which vanishes because of the cobordism invariance of the signature. 
Therefore using \eqref{indin} we get the following formula for relative Kervaire semi-Characteristic  
\begin{equation}\label{finest}
\kappa(M,\partial M)=\sum_{\gamma\subset M^\circ}\ind_2(\gamma)
\end{equation}
and this complete the proof of the theorem \ref{aslio}.

\section{Application}
\subsection{Cut-and-past invariance of Kervaire semi-Characteristic}
In this section we use the main theorem \ref{aslio} to prove the cut-and-past invariance of the 
Kervaire semi-characteristic. Let $M$ be a closed manifold without 
boundary and let $N\subset M$ 
be a sub-manifold of codimension one that partitions $M$ into two compact manifolds $M_1$ and $M_2$ 
with common boundary $\partial M_1=\partial M_2=N$. Given a diffeomorphism  $\phi: N\to N$ the closed manifold 
$M_1\sqcup_\phi M_2$ is defined by identifying points $x\in \partial M_1$ and $\phi(x)\in\partial M_2$. 
It is clear that $M_1\sqcup_\phi M_2$  is produced by cutting $M$ along $N$ and then pasting through $\phi$, 
hence the nomination. 

\begin{thm}\label{theorem2}
Let $M$ be a $(4k+1)$-dimensional closed manifold. 
Kervaire semi-characteristic is invariant with respect to cut-and-past operation 
$\kappa(M)=\kappa(\partial M_1\sqcup_\phi\partial M_2)$ provided that $\chi(N)=0$. 
\end{thm}  
\begin{pf}
Since $\chi(N)=0$ and considering the dimensions, there is a nowhere vanishing vector field $V$ on $N$ 
that has a nowhere vanishing extension to $M_1$. We denote this extension by $V_1$.   
Let $E_1$ be the bundle over $M_1$ whose fiber at $x$ is orthogonal to $V_1(x)$. Let $X_1$ be a 
transversal section of $E$ which is tangent to $N$. The zero set of $X$ consists some 
circles $\gamma$ in $M_1$ 
which are either disjoint from, or included in $N$. Moreover each of these circle carries a real line bundle  
$\mathcal L$. Now consider the push-forward vector field $V_2=\phi_*(V_1)$ and extent it 
to a nowhere vanishing 
vector field $V_2$ on $M_2$. Let $E_2$ be the vector bundle on $M_2$ whose fiber are orthogonal 
to $V_2$. Let $X_2$ 
be a transversal section of $E_2$ such that on the boundary $N$ it satisfies $X_2=\phi_*(X_1)$. 
These geometric 
data fit together to build up corresponding geometric data on the whole manifold 
$\partial M_1\sqcup_\phi\partial M_2$. This way we get vector fields $V$ and $X$ on, and vector 
bundle $E$ over  $\partial M_1\sqcup_\Phi\partial M_2$. Applying  the Zhang counting formula \eqref{zhangasli} 
to these data we get the following relation 
\[\kappa(\partial M_1\sqcup_\phi\partial M_2)=\sum_{\gamma\in M_1^\circ}\ind_2(\gamma)
+\sum_{\gamma\in M_2^\circ}\ind_2(\gamma)+\sum_{\gamma\in N}\ind_2(\gamma).\]
Using the theorem \ref{aslio}, the first and second terms in the right side of the above 
expression are respectively equal to relative Kervaire semi-characteristics $\kappa(M_1,\partial M_1)$ 
and $\kappa(M_2,\partial M_2)$. Concerning the last term, it vanishes because of the cobordism invariance of signature. Therefore  
\begin{equation}
\kappa(\partial M_1\sqcup_\phi\partial M_2)=\kappa(M_1,\partial M_1)+\kappa(M_2,\partial M_2). 
\end{equation}
The right side of this relation has no reference to the pasting diffeomorphism $\phi$ so 
the left side is indeed independent of $\phi$.  
\end{pf}

\end{document}